\documentclass[11pt]{article}
\usepackage{amsfonts}
\usepackage{setspace}
\usepackage{amssymb}
\usepackage{amsmath}
\usepackage{eufrak, makeidx}
\usepackage{mathrsfs}
\usepackage{amssymb, enumerate, amsthm, amsfonts, amsmath, mathrsfs}
\usepackage{natbib}
\newtheorem{thm}{Theorem}

\newtheorem{dfn}[thm]{Definition}

\usepackage{rotating}

\makeindex
\begin{document}
	\pagenumbering{arabic} 
	\setcounter{page}{1}
	\large 
	\begin{center}
		\textbf{A direct construction of the Wiener measure on $\textbf{C}[0, \infty)$}\\[0.6 cm]
		R. P. Pakshirajan and M. Sreehari*\\[0.5 cm] \end{center}
	\begin{footnotesize}
		\makeatletter{\renewcommand*{\@makefnmark}{}
			\footnotetext{R. P. Pakshirajan, 227, 18th Main, 6th Block, Koramangala, Bengaluru -560095, India.\\
				*Corresponding author: M. Sreehari, 6-B, Vrundavan Park, New Sama Road, Vadodara, 390024, India.\\	
				
				E-\textit{mail addresses}: vainatheyarajan@yahoo.in (R. P. Pakshirajan);  msreehari03@yahoo.co.uk (M. Sreehari)  }}
		\end {footnotesize}
		
		\noindent	\textbf{Abstract.}\\
		Our construction of the Wiener measure on $\textbf{C}=\textbf{C}[0, \infty)$ consists in first defining a  set function $\varphi$\ on the class of all compact sets based on certain $n$-dimensional normal distributions, $n = 1,\ 2,\ldots$\ using the structural relation at (\ref{E1.2}) below. This structural relation, discovered by the first author, is recorded in his book (2013) on page 130. We then define a measure $\mu$ on the  Borel $\sigma$-field of subsets of $\textbf{C}$ which is the Wiener measure.  This is done via a similar construction of the Wiener measure on  $\textbf{C}_a=\textbf{C}[0, a)$  where $a > 0$ is an arbitrary real number.\\
		The traditional way is to first construct the Brownian Motion process (BMP) and then, by proving it is a measurable mapping into $(\textbf{C},\ \mathscr{C}_\infty)$, call the measure induced by the BMP on $\textbf{C}$\ the Wiener measure. In the present paper, we define the Wiener measure directly. 
	\section{INTODUCTION AND PRELIMINARY RESULTS.}
	Construction of the Wiener process is discussed by many authors and the discussion invariably starts with the construction first of the Brownian Motion Process (BMP) on a probability space. The BMP is studied for its properties and then is proved to be a measurable mapping into $\textbf{C}[0,\ 1]$ space 
endowed with the uniform metric and the resulting Borel $\sigma$-field $\mathfrak{C}_1$. Call the measure induced by the BMP on $\textbf{C}[0, 1]$\ the Wiener measure. 
We refer to chapter 2 in  Karatzas and Shreve (1988).\\ The aim of the present paper is to reverse the procedure and construct the Wiener measure  directly  using elementary measure  theory  and the structural relation given at (\ref{E1.2}) below.  In Pakshirajan and Sreehari, 2020 the authors presented the construction of the Wiener measure on \textbf{C}[0, 1].\\
	Let $a > 0$\ be arbitrary but fixed. Let $\textbf{C}_a = \textbf{C}[0,\ a]$\ denote the space of real valued continuous functions defined on $[0,\ a]$, all vanishing at $0$\ and endowed with
the norm $\|x\|_a = \sup\limits_{0 \le t \le a}|x(t)|,\ x \in \textbf{C}_a$. Let $\rho_a(x,\ y) $\ denote the associated metric. Define another norm : $\|x\|^*_a = \sup\limits_{0 \le s,\ t \le a}|x(t) - x(s)|$. The associated metric will be denoted by $\rho^*_a$. Since $\|x\|_a \le \| x \|^*_a\le 2\| x \|_a$, the two norms induce the same topology in $\textbf{C}_a$\ and determine the same Borel $\sigma$-field.\\
\noindent Let $T_n = \{\frac{{a}k}{2^{n}},\ k = 0,\ 1,\ 2,\ldots,\ 2^{n}\}; T = \bigcup\limits_{n = 1}^{\infty}T_n$\ and note that $T$\ is a countable dense subset of the interval $[0,\ a]$. For $x \in \textbf{C}_a$, define
\begin{equation} \label{E2.3}
	\wp_nx = \Big(x(\frac{{a}}{2^n}),\ x(\frac{2{a}}{2^n}) - x(\frac{{a}}{2^n}),\ldots, x(\frac{{a}2^{n}}{2^n}) - x(\frac{{a}(2^{n} - 1)}{2^n})\Big).
\end{equation}
This maps $\textbf{C}_a$\ into $R^{2^{n}}$. Assume $R^{2^{n}}$\ is endowed with the usual metric and denote the resulting Borel $\sigma$-field by $\mathfrak{R}^{2^{n}}$. We note $\wp_n$\ is a continuous map and hence is $\mathfrak{C}_{a} $\ measurable. We prescribe the distribution of the vector variable $\wp_n$\ to be the multivariate normal distribution with independent components, each component with zero mean and variance $\frac{{a}}{2^n}$. i.e., it is the joint distribution of $(\sqrt{\frac{{a}}{2^n}}\xi_k,\ 1 \le k \le 2^n)$\ where the $\xi_k$s are independent standard normal variables.\\ 
Denote by $\nu_n$\ the measure on $\mathfrak{R}^{2^n}$\ by this distribution.\\ Let $\alpha_{n}$\ denote the measure generated on the sub $\sigma$-field
$\wp_n^{-1}(\mathfrak{R}^{2^n})$\ by the mapping $\wp_{n}$. 
All sets considered below are members of $\mathfrak{C}_a$.\\
\noindent Let $K \subset \textbf{C}_a$\ be compact. Then the following structural relation holds:
(ref. pp 130-131 in Pakshirajan, 2013.)\
\begin{equation}\label{E1.2}
	K = \cap_{n = 1}^{\infty}\wp_n^{-1}\wp_nK.
\end{equation}
To make for seamless reading we present here a proof of (\ref{E1.2}).\\
\noindent That $K \subset \cap_{n = 1}^{\infty}\wp_n^{-1}\wp_nK$\ is obvious. Now to establish the reverse inclusion, let $x$\ be an arbitrary member of the right side. Hence for every $n$, $x \in \wp_n^{-1}\wp_nK$. There exists therefore $y_n \in K$\ such that $\wp_nx = \wp_ny_n$. Since $K$\ is compact, sequence $(y_n)$\ contains a convergent subsequence, say, $(y_m)$\ converging to say $y_0 \in K$\ in the metric $\rho_{\mathfrak{a}}$. This implies $y_m(t) \rightarrow y_0(t)$\ for all $t \in [0, a]$. 
Fix $r$\ and let $1 \le j \le 2^{r}$.  Let $m > r$. The relation $\wp_mx = \wp_my_m$, is equivalent to the relation $\wp_mx = \wp_my_m$\ in the sense that given $\wp_mx \in R^{2^m}$\ the point $\wp_mx$\ is uniquely determined and conversely through a linear transformation. Here $\wp_mx = \Big(x(\frac{{a}}{2^m}),\ x(\frac{2{a}}{2^m}),\ ...,\ x(\frac{{a}2^m}{2^m})\Big)$. We get
$x(\frac{{a}j}{2^r}) = y_m(\frac{{a}j}{2^r})$. Take limit as $m \rightarrow \infty$, and get $x(\frac{{a}j}{2^r}) = y_0(\frac{{a}j}{2^r})$. Thus for every $u \in T,\ x(u) = y_0(u)$. Since $T$\ is dense in $[0,\ {a}]$\ and since $x,\ y_0$\ are continuous functions, it follows that $x(t) = y(t)$\ for all $t \in [0,\ {a}]$. Thus $x \in K$\ and the proof is complete.\\ 
Note that this inclusion is true for any set $K$ and not only for compact sets.
\begin{thm}\label{T1.1}
	For any $A \in \mathfrak{C}_a,\ \alpha_{n}(\wp_{n}^{-1}\wp_{n}(A)),\ n = 1,\ 2,\ldots.$\ is a
	monotonic decreasing sequence of numbers.
\end{thm}
\noindent	\textbf{Proof.}\\
\begin{equation*}
	\begin {split}
	\alpha_{n + 1}\big(\wp_{n + 1}^{-1}\wp_{n + 1}(A)\big) = \int_{\wp_{n + 1}(A)}\;d\nu_{n + 1} \le \int\limits_{\wp_n(A) \times R}\;d\nu_n\;d\beta_{n + 1}\\
	\le \int\limits_{\wp_n(A)}\;d\nu_n \le \alpha_n(\wp_n^{-1}\wp_n(A))
\end{split}
 \end{equation*} 
where $\beta_{n +1}$ is the distribution function of a normal variable and $\nu_n=\alpha_n \wp_n^{-1}.$ \\
\noindent Define set function $\varphi$\ on the compact sets $K$\ of $\textbf{C}_a$ :
\begin{equation}  \label{E1.3}
	\varphi(K) = \lim_{n \rightarrow \infty}\alpha_{n}\big(\wp_n^{-1}\wp_n(K)\big).
\end{equation} 
Note \begin{equation}\label{E1.4}
	\varphi(K) \le 1;\ \ \varphi(\emptyset) = 0.
\end{equation}
\begin{thm} \label{T1.2}
	Let $K_1,\ K_2$\ be compact sets such that $\varphi(K_1 \cap K_2) = 0$. Then $\varphi(K_1 \cup K_2)= \varphi(K_ 1) + \varphi(K_2)$.
	
\end{thm}
\noindent	\textbf{Proof.}\\
$\wp_n(K_1 \cup K_2) = \wp_n(K_1) \cup \wp_n(K_2))$\\
$ \wp_n^{-1}\wp_n(K_1 \cup K_2) = \wp_n^{-1}\wp_n(K_1) \cup \wp_n^{-1}\wp_n(K_2)$.\\
Hence
\begin{equation*}
	\begin{split}	
		\alpha_n\big(\wp_n^{-1}\wp_n(K_1 \cup K_2)\big) = \alpha_n\big(\wp_n^{-1}\wp_n(K_1)\big) + \alpha_n\big(\wp_n^{-1}\wp_n( K_2)\big) \\- \alpha_n\big(\wp_n^{-1}\wp_n(K_1) \cap \wp_n^{-1}\wp_n(K_2)\big)
\end{split} \end{equation*}
since $\alpha_n$\ is a measure. Now, $\wp_n^{-1}\wp_n(K_1) \cap \wp_n^{-1}\wp_n(K_2) = \wp_n^{-1}\wp_n(K_1 \cap K_2)$.
Since $K_1 \cap K_2$\ is a compact set and since $\varphi(K_1 \cap K_2) = 0$, $\alpha_n\big(\wp_n^{-1}\wp_n(K_1 \cap K_2)) < \varepsilon$\ for all $n$\ large. Taking limits as $n \rightarrow \infty$\ and then as $\varepsilon \rightarrow 0$\ in the inequalites
\begin{equation*}
	\begin{split}
		\alpha_n\big(\wp_n^{-1}\wp_n(K_1)\big) + \alpha_n\big(\wp_n^{-1}\wp_n( K_2)\big)  - \varepsilon \le \alpha_n\big(\wp_n^{-1}\wp_n(K_1 \cup K_2)\big) \\
		\le  \alpha_n\big(\wp_n^{-1}\wp_n(K_1)\big) + \alpha_n\big(\wp_n^{-1}\wp_n( K_2)\big)	
\end{split}\end{equation*} 
we complete the proof of the claim.\\
\noindent	\textbf{Remark.}\label{R1.1}	
	We have the following observations from the earlier discussion:\\
	a) $\varphi$ is finitely additive on the collection of compact sets.\\
	b) $0 \le \varphi (K) \leq 1$ for all compact sets $K$.\\
	c) If $K_1, K_2$ are compact sets and $K_1 \subset K_2$ then from (\ref{E1.3}),  $\varphi(K_1) \leq \varphi(K_2)$.\\

\begin{dfn} We call a set, in a topological space, a boundary set if it is a closed set with a null interior. The boundary of  a set $A$ (i.e.,   $\overline{A} \sim\textit{Int}\;A$ ) will be denoted by  $\partial A$. \\
	We note that the boundary of a set is a boundary set.
\end{dfn}
\begin{thm}\label{T1.3}
	If $K_1,\ K_2$\ are compact subsets with $K_1 \subset K_2,\ \text{and}\ \varphi(\partial{K_1}) = 0$, then $\varphi(K_2 \cap \bar{K_1^{\prime}}) = \varphi(K_2) - \varphi(K_1)$. 
\end{thm}
\noindent \textbf{Proof.}  $\bar{K_1^{\prime}} \cap K_2$\ is a compact set. $K_1 \cap \{\bar{K_1^{\prime}} \cap K_2\} = \partial K_1$. Since $\varphi(\partial K_1) = 0$, Theorem \ref{T1.2} applies and we get $\varphi(K_2) = \varphi\big(K_1 \cup \{\bar{K_1^{\prime}} \cap K_2\}\big) = \varphi(K_1) + \varphi\big(\bar{K_1^{\prime}} \cap K_2\big)$, as was to be proved. \\
\noindent We now discuss some limiting properties of $\varphi(K_n)$.
\begin{thm} \label{T1.4}
	(i) Let $K_n, n \geq 1$, be compact sets such that $K_n \downarrow K$. Then $\varphi(K_n) \downarrow \varphi(K).$\\
	(ii) Suppose $K, K_n, n \ge 1$ are compact subsets with $K_n \uparrow K\ \text{and}\; \varphi(\partial K_n) = 0$.  Then $\varphi(K_n) \uparrow \varphi(K)$.\\
	(iii) Let $K,\ K_m,\ m \ge 1$\ be compact sets with $K_m \uparrow K$\ and $\varphi(K_m) = 0$. Then $\varphi(K) = 0$.	
\end{thm}
\noindent \textbf{Proof.}
(i) Since the sequence $(\varphi(K_n))$\ is monotonic decreasing, it is enough to show that, given $\varepsilon > 0$, there exists $K_N$\ such that $\varphi(K_N) < \varphi(K) + \varepsilon$.\\
We note $K$\ is compact. Hence given $\varepsilon > 0$, we can find $r \ge 1$\ such that
\begin{equation}\label{E1.5}
	\varphi(K) > \alpha_{\ell}\big(\wp_{\ell}^{-1}\wp_{\ell}(K)\big) - \varepsilon
\end{equation}
for all $\ell \ge r.$ Since $K_n \downarrow K$, for all $\ell \geq 1$ we have $\wp_{\ell}^{-1}\wp_{\ell}(K_n) \downarrow \wp_{\ell}^{-1}\wp_{\ell}(K)$.\\
For fixed $\ell$ we then have, as $n\rightarrow \infty$
$\alpha_{\ell}\big(\wp_{\ell}^{-1}\wp_{\ell}(K_n)\big) \downarrow \alpha_{\ell}\big(\wp^{-1}_{\ell}\wp_{\ell}(K)\big)$.
Take $\ell = r$. We can find $N = N(r)$\ large such that
$\alpha_r\big(\wp^{-1}_r \wp_r(K_N)\big) < \alpha_r\big(\wp^{-1}_r\wp_r(K)\big) + \varepsilon.$
This, together with (\ref{E1.5}), yields
\begin{equation*}
	\varphi(K) + \varepsilon > \alpha_r\big(\wp^{-1}_r \wp_r(K_N)\big) - \varepsilon > \varphi(K_N) - \varepsilon.
\end{equation*}
Since $\varepsilon > 0$\ is arbitrary, the claim follows.\\	
(ii) Define $E_n = K \cap \bar{K_n^{\prime}}$\ and note by Theorem  \ref{T1.3} that $\varphi(E_n) = \varphi(K) - \varphi(K_n)$. Now, the $E_n$s are compact sets and $E_n \downarrow \emptyset$. Hence by part (i) above, $\varphi(E_n) \downarrow 0$. i.e., 
$\varphi(K_n) \rightarrow \varphi(K)$, as was to be proved.\\
(iii) Claim immediate from part(ii) above.
	\section {THE WIENER MEASURE on $\textbf{C}_a$.}
In this Section we introduce a new set function in terms of $\varphi$ on the Borel $\sigma$-field $\mathfrak{C}_a$ of subsets of $\textbf{C}_a$ and study its properties to show that it is indeed the Wiener measure.\\
For  arbitrary measurable sets  $A \in \mathfrak{C}_a$\, define
\begin{equation}	\label{E2.1}
	\mu(A) = \sup\limits_{K \subset A,\ K\ \text{compact}}\varphi(K).
\end{equation}
\noindent At the outset we observe that for compact sets $K$, $\mu(K)=\varphi(K)$ and hence all the properties noted in the previous Section for $\varphi$ also hold for $\mu$. Further the definition 
implies (i) that if $A \subset B, A, B \in \mathfrak{C}_a\ \text{then}\ \mu(A) \le \mu(B)$\ and (ii) that there exists an increasing sequence ($K_n$) of compact sets, $K_n \subset A$\ such that $\mu(A) = \lim\limits_{n \rightarrow \infty}\mu(K_n)$.\\
The sets $K_n$ can be chosen to be monotonic increasing.\\
	\textbf{Remark.}
	This does not mean that $K_n \uparrow A$. i.e., $\bigcup\limits_{n = 1}^{\infty}K_n$\ can be a proper subset of $A$. 
	To see this, take $v \in \textbf{C}, \|v\|=1.$ Let $K_n=\{\lambda v,\; 0 \leq \lambda \leq 1-\frac{1}{n}\}$ and  $A=\{\lambda v, \; 0\leq \lambda \leq 1 \}.$  However, if  $K_n=\{\lambda v,\; 0 \leq \lambda \leq 1-\frac{1}{n}\}\cup \{v\}$ then both $K_n$ and $A$ are compact and  $K_n\uparrow A$.\\

\noindent We next discuss further properties of $\mu$ that enable us to claim that $\mu$ is indeed a probability measure.
\begin{thm}\label{T2.1}
	(i) If $A \subset B$, then $\mu(A) \le \mu(B)$.\\
	(ii)If $A \cap B = \emptyset$, then $\mu(A \cup B) = \mu(A) + \mu(B)$.
\end{thm}
\noindent \textbf{Proof.}
(i) immediate from the definition of $\mu$\ at (\ref{E2.1})\\
(ii) Let $E \subset A,\ F \subset B$\ be compact sets such that for a given $\varepsilon > 0$ \; $\mu(E) > \mu(A) -\varepsilon$ and $\mu(F) > \mu(B) - \varepsilon $. We have, from Theorem \ref{T1.2}\\
$\mu(E) + \mu(F) = \mu(E \cup F) \le \mu(A \cup B)$ since $\varepsilon$ is arbitrary.\\
Thus $\mu(A \cup B) \ge \mu(A) + \mu(B)$. It remains to be shown that $\mu(A \cup B) \le \mu(A) + \mu(B).$\\
Given $\varepsilon > 0$, we can find a compact set $K,\;\ K \subset A \cup B$\ such that\\ 
$\mu(A \cup\ B) - \varepsilon < \mu(K)$.\\
Case 1. The distance  $d(A,\ B) = q > 0$.\\ Consider an arbitrary sequence $(x_n)\ \text{in}\ K \cap A$. Since it is a sequence in $K$, it contains a convergent subsequence, converging to, say, $x_0$. This $x_0$\ has to be in $K \cap A$\ or in $K \cap B$. Since the sequence lies in $K \cap A$\ and since $d(K \cap A,\ K \cap B) \ge q > 0$, we conclude $x_0 \in K \cap A$. Thus we see that every sequence in $K \cap A$\   
contains a convergent subsequence converging to a point in $K \cap A$. This means $K \cap A$\ is a compact set. Similarly, $K \cap B$\ is a compact set. Summarising, we conclude that every compact subset of $A \cup B$\ is the union of a compact subset $E$\ of $A$\ and a compact subset $F$\ of $B$. We get $\mu(A \cup\ B) - \varepsilon < \mu(K) = \mu(E \cup F) = \mu(E) + \mu(F) \le \mu(A) + \mu(B)$.  That $\mu(A \cup B) \le \mu(A) + \mu(B)$\ is now immediate.\\ 
Case 2. $d(A,\ B) = 0$.\\ This case assumption implies that $Q = \bar{A} \cap \bar{B} \neq \emptyset$. Again in this case one or both the sets $K \cap A,\ K \cap B$\ can fail to be compact. Since the other case admits to being similarly argued, let us assume that neither of the two sets is compact. 
$K \subset A \cup B$\ can not be compact if any convergent sequence in it converges to a point outside $K$. i.e., if convergent sequences in $K \cap A$\ or in $K \cap B$\ converge to points outside these sets. Thus $K$\  can be a compact subset only if $E = K \cap A$\ and $F = K \cap B$\ are compact. And the arguments and the conclusion in case 1 hold.\\ With this the proof is complete.\\
	\textbf{Remark.}	\label{R2.2}
	Immediate consequences of Theorem \ref{T2.1} are :\\
	a) If $A_k,\ 1 \le k \le n$\ is any collection of $n$\ events, then $\mu(\bigcup\limits_ {k = 1}^nA_k) \le \sum\limits_{k = 1}^n\mu(A_k)$ and if the events $A_n$ are mutually exclusive equality holds.\\
	(b) If $A \subset B$, then $\mu(B \sim A) = \mu(B) - \mu(A)$.\\

Our next result shows that $\mu$ is monotone.
\begin{thm} \label {T2.2}
	(i) If $A_n \downarrow A$, then $\mu(A_n) \downarrow \mu(A)$.\\
	(ii)If $A_n \uparrow A$, then $\mu(A_n) \uparrow \mu(A)$.
\end{thm}
\noindent\textbf{Proof.}\\
(i) By Theorem \ref{T2.1} we note that the hypothesis implies $B_n \downarrow \emptyset$\ where $B_n = A_n \cap A^{\prime}$. 
	We refer to Remark \ref{R2.2}(a) and claim that it is enough to show that $\mu(B_n) \rightarrow 0$.\\
Find compact sets $K_n \subset B_n$\ such that $\mu(B_n) - \mu(K_n) < \frac{\varepsilon}{2^n}$. Define $Q_n = \cap_{j = 1}^nK_j$, Note that $Q_n \subset B_n$, that $Q_n$\ is a compact set and that $Q_n \downarrow \emptyset$, By Theorem \ref{T1.4} and the fact that $\mu(\phi)=0$ by (\ref{E1.4}), it then follows that $\mu(Q_n) \rightarrow 0$. Further by Remark \ref{R2.2}(a)
\begin{equation*}
	\begin{split}
		\mu(B_n) - \mu(Q_n) = \mu(B_n \cap Q_n^{\prime})\  = \mu\big(B_n \cap \{\bigcup\limits_{j = 1}^nK_j^{\prime}\}\big) = \mu\big(\bigcup\limits_{j = 1}^n(B_n \cap K_j^{\prime})\big)\\
		\le \sum\limits_{j = 1}^n\mu(B_n \cap K_j^{\prime}) \le \sum\limits_{j = 1}^n\mu(B_j \cap K_j^{\prime})\le \sum\limits_{j = 1}^n\{\mu(B_j) - \mu(K_j)\}\le \sum\limits_{j = 1}^n\frac{\varepsilon}{2^j} < \varepsilon
	\end{split}
\end{equation*}
for all $n$. Here we used the Remark \ref{R2.2} and the fact that $B_n$  is $\downarrow$. Collecting the results, we conclude $\mu(B_n) \rightarrow 0$, thus completing the proof of this part.\\
(ii) That $\mu(A_n) \uparrow$\ is true follows from Theorem \ref{T2.1}(i). Since $A \cap A_n^{\prime} \downarrow \emptyset$, part (i) applies and we have $\mu(A_n \cap A^{\prime}) \rightarrow 0$. Now by Remark \ref{R2.2}(b), this gives $\mu(A_n) \rightarrow \mu(A)$.
\begin{thm} \label{T2.3}
	$	\mu$ defined at (\ref{E2.1}) is a probability measure.
\end{thm}
\noindent\textbf{Proof.}\\
Let $A_n \in \mathfrak{C}_a,\ n \ge 1$\ be a sequence of mutually exclusive events. Let $A = \bigcup_{n = 1}^{\infty}A_n = \bigcup_{n = 1}^{\infty}B_n$\ where $B_n = \bigcup_{k = 1}^nA_k$. Since $B_n \uparrow A$, Remark \ref{R2.2}(a) applies and then we have 
\begin{equation*}
	\mu(A) = \lim_{n \rightarrow \infty}\mu(B_n) = \lim_{n \rightarrow \infty}\sum_{k = 1}^n\mu(A_k)= \sum\limits_{k = 1}^{\infty}\mu(A_k).
\end{equation*}
i.e., $\mu$\ on $\mathfrak{C}_a$\ is countably additive. Since $\mu (A) \geq 0$ for $A\in \mathfrak{C} $  it follows that $\mu$ is a probability measure if we show that $\mu(\textbf{C}_a)=1$.\\
Let $T = \cup_{n = 1}^{\infty}T_n,\ T_n = \{t_k,\ 1 \le k \le 2^n\}$\ where $t _k = t_{k,n} = \frac{ak}{2^n}$\ and note that $T$\ is a countable dense subset of the  interval $[0,\ a]$. Let $S_m = \{x:\ x \in \textbf{C}_a,\ \Arrowvert x \Arrowvert \le m\}$. We note $S_m = \{x:\ x \in \textbf{C}_a,\ \sup\limits_ {t \in T}|x(t)| \le m\} = \cap_{n = 1}^{\infty}B_{n,m} = \lim\limits_{n \rightarrow \infty}B_{n,m}$\ where $B_{n,m} = \{x:\ x \in \textbf{C}_a,\ \sup\limits_{t \in T_n}|x(t)| \le m\} = 
\cap_{t \in T_n}\{x:\ x \in \textbf{C}_a,\ |x(t)| \le m\}$.\\
We note that $S_m \uparrow \textbf{C}_a$. Recall that, given $\varepsilon > 0$, we can find $A_m \subset S_m,\ A_m$\ compact such that $\mu(S_m) - \mu(A_m) < \varepsilon$. Write\\ $\wp_{T_n}x = \Big(x((\frac{a}{2^n}),\ x(\frac{2a}{2^n}) - x(\frac{a}{2^n}),\ x(\frac{3a}{2^n}) - x(\frac{2a}{2^n}),\ ...,\ x(\frac{a2^n}{2^n}) - x(\frac{a(2^n - 1)}{2^n})\Big).$\\
If $K \subset \textbf{C}_a$\ is a compact set, then arguing as in the proof of Theorem  \ref{T1.1} we get  $\wp^{-1}_{T_n}\wp_{T_n}K \downarrow K$.
Hence given $\varepsilon > 0$, we can find $N$\ such that for all $n \ge N$,\\
$\mu(A_m) + \varepsilon > \alpha_{T_n}\wp^{-1}_{T_n}\big(\wp_{T_n}S_m\big) = P\big(\frac{a^{1/2}}{2^{n/2}}\max\limits_{1 \le j \le 2^{n}}|\xi_j| \le m\big)$ \\
where the $\xi$s are independent standard normal variables. Hence \\
$\mu(A_m) + \varepsilon > \big(P(|\xi| \le \frac{2^{n/2}\;m}{a^{1/2}})\big)^{2^n} = \big(1 - P(|\xi| > \frac{2^{n/2}\;m}{a^{1/2}})\big)^{2^n} \ge \big(1 -  \frac{a\mathbb{E}|\xi|^2}{m^2\;2^n}\big)^{2^n} $\\
leading to $\mu(A_m) + \varepsilon \ge e^{- (a/m^2)}.$\\This implies by Theorem \ref{T2.2}(ii),\\ $\mu(\textbf{C}_a) = \lim\limits_ {m \rightarrow \infty}\mu(S_m) \ge \lim\limits_ {m \rightarrow \infty}\mu(A_m) \ge \lim\limits_{m \rightarrow \infty}e^{- a/m^2} - \varepsilon \ge 1 - \varepsilon$.\\
Since $\varepsilon > 0$ is arbitrary we get $\mu(\textbf{C}_a)=1.$\\
\textbf{	Alternate proof for $\mu(\textbf{C}_a) = 1$ which will be useful in Section 3.}\\
Since $\mu(K) \le 1$\ for all compact sets, as noted in (\ref{E1.4}) and since \\$\mu(\textbf{C}_ {{a}}) = \sup\limits_{K\subset \textbf{C}_{\mathfrak{a}},\ K\ \text{compact}}\mu(K)$, it follows that $\mu(\textbf{C}_{{a}}) \le 1$. So the proof will be complete if we show that $\mu(\textbf{C}_{{a}}) \ge 1$. This we proceed to show.\\ 
	Let
\begin{equation*}
	H^0_\alpha=\{x: \sup_{0 \le s, t \le 1;\;s \ne t} \frac{|x(t)-x(s)|}{|t-s|^\alpha} < \infty\}  
\end{equation*}
and 
\begin{equation*}
	H^0_{\alpha, a}=\{x: \sup_{0 \le s, t \le a;\;s \ne t} \frac{|x(t)-x(s)|}{|t-s|^\alpha} < \infty\} .
\end{equation*}
Let $\delta > 0$ be arbitrary. Then for $|s-t| < \delta$ and $ x\in H^0_\alpha$  without loss of generality we have
\begin{equation*}
	|x(s)=-x(t)| < |t-s|^\alpha < \delta ^\alpha.
\end{equation*}
Then we have the following Theorem which in turn implies $\mu(\bf{C}_a ) =1.$
	\begin{thm}
	For $0 < \alpha < 1, \;\;\mu(H^0_\alpha) =1.$
\end{thm}
\noindent \textbf{Proof.}\\
Take $ n$ large so that $\frac{a}{2^n} < \delta $. Then $|x(\dfrac{a(r+1)}{2^n}) -  x(\dfrac{a{r}}{2^n}) | \le (\frac{a}{2^n})^\alpha$ for $r=0, 1, \ldots, n-1.$\\
Since the $\mu$ measure of every compact subset of $\bf{C}_a  \le 1$ it follows that the $\mu$ measure any measurable subset of $\bf{C}_a  \le 1$  as well. Now since $$	H^0_{\alpha, a}=\{x: \sup_{0 \le s, t \le a;\;s \ne t} \frac{|x(t)-x(s)|}{|t-s|^\alpha} < \infty\}$$ is a measurable subset of  $\bf{C}_a$ it follows that $\mu(	H^0_{\alpha, a}) \le 1.$ \\
Set \begin{equation*}
	\mathbb{S}_{\alpha, a}(\lambda)=\{ x: \sup_{0< s, t \le a;\; 0 < |s-t| < 1}\frac{|x(t)-x(s)|}{|s-t|^\alpha} \le \lambda\}.
\end{equation*}
Then $\mathbb{S}_{\alpha, a}(\lambda)$ is a compact subset of $\bf{C}_a$ (see Appendix) and $\mu ((\mathbb{S}_{\alpha, a}(\lambda))^\prime) < \varepsilon$ for all large $\lambda$ depending on $\varepsilon$.	i.e., $\lim_{n \rightarrow \infty} \nu_n(\wp_n(\mathbb{S}_{\alpha, a}(\lambda))^\prime)) < \varepsilon.$\\ Let 
\begin{equation*}
	A_n= \{x: x \in H^0_{\alpha, a}; \max_{0 \le r \le 2^n-1}|x(a(r+1))-x(ar)|\le \frac{\lambda a^\alpha}{2^{n \alpha}}\}.
\end{equation*}
Note that for $x \in A_n$ 
$$\wp_n(x)= \bigg(x(\frac{a}{2^n}), x(\frac{2a}{2^n})-x(\frac{a}{2^n}), \ldots, x(\frac{a2^n}{2^n})-x(\frac{a(2^n -1)}{2^n})\bigg) $$
We then have 
\begin{equation*}
	\begin{split}
		\mu(\mathbb{S}_{\alpha, a}(\lambda))=\lim_{n \rightarrow \infty} \nu_n(\wp_n(\mathbb{S}_{\alpha, a}(\lambda)))\ge \lim_{n \rightarrow \infty} \nu_n(\wp_n(\mathbb{S}_{\alpha, a}(\lambda)\cap A_n))\\
		\ge	\lim _{n \rightarrow \infty}[\nu_n(\wp_nA_n)-\nu_n(\wp_n(A_n\cap(\mathbb{S}_{\alpha, a}(\lambda))^\prime)]\\
		\ge	\lim _{n \rightarrow \infty}[\nu_n(\wp_nA_n)-\nu_n(\wp_n((\mathbb{S}_{\alpha, a}(\lambda))^\prime)]\\
		\ge	\lim _{n \rightarrow \infty}[\nu_n(\wp_nA_n)-\varepsilon]\\
		\ge \lim_{n \rightarrow \infty} P( \max_{0 \le r \le 2^n-1}|x(a(r+1))-x(ar)| \le \frac{\lambda a^\alpha}{2^{n \alpha}})-\varepsilon\\
		\ge\lim_{n \rightarrow \infty} P(\max_{0 \le r \le 2^n-1}|\xi_r| \le \frac{\lambda 2^{n(1-2\alpha)/2}}{a^{(1-2\alpha)}/2})-\varepsilon\\
		=\lim_{n \rightarrow \infty}[P(|\xi|\le\frac{\lambda 2^{n(1-2\alpha)/2}}{a^{(1-2\alpha)}/2})]^{2^n}-\varepsilon\\
		=\lim_{n \rightarrow \infty}[1-P(|\xi|>\frac{\lambda 2^{n(1-2\alpha)/2}}{a^{(1-2\alpha)/2}})]^{2^n}-\varepsilon\\
		\ge \lim_{n \rightarrow \infty}[1- \frac{a}{2^n}\frac{E|\xi|^{2/(1-2\alpha)}}{\lambda^{2/(1-2\alpha)}}]^{2^n}-\varepsilon\\
		\ge e^{-\psi(\lambda)}- \varepsilon	
	\end{split}				
\end{equation*}
by Tchebichev's inequality 	where $\xi, \xi_k$ are independent standard normal rvs and $\psi(\lambda)=\frac{a E|\xi|^{2/(1-2\alpha)}}{\lambda^{2/(1-2\alpha)}} \rightarrow 0$ as $\lambda \rightarrow \infty$. Since $\varepsilon$ is arbitrary From the above result we get  $	\mu(\mathbb{S}_{\alpha, a}(\lambda))\ge 1$ and  hence $\mu(H^0_\alpha) =1$. This completes an alternate proof of Theorem 11.\\
	\textbf{Remark.} \label{R2.3}
	(i) From the construction of $\mu_{{a}}$, it is clear that if $\nu$\ is a probability measure on $\mathfrak{C}_{{a}}$\ and if its finite dimensional distributions (i.e., the distributions of the vector variables ($\pi_{t_1},\ \pi_{t_2},\ldots,\ \pi_{t_k}$), for every choice of $k$\ and every choice of ($t_1,\ t_2,\ldots,\ t_k$) are the same as the corresponding ones of $\mu_{{a}}$\ then $\nu \equiv \mu_{{a}}$. It follows
	now that $\mu$ is the Wiener probability measure.\\
	(ii) The co-ordinate process $\{\pi_t,\ t \ge 0\}$\ is known as the Brownian motion process. \\

	\section{Constructing the Wiener measure on $\textbf{C}_{\infty}.$ }
Let $\textbf{C}_{\infty}$ be the space of all the real continuous functions defined on $[0, \infty)$, all vanishing at $0$ endowed with the metric  $$d(f, g)= \sum_{n=1}^\infty \frac{1}{2^n} \frac{\sup_{0\le  t \le  n}|f(t)-g(t)|}{1+\sup_{0\le  t \le n}|f(t)-g(t)|}$$
Further on $\mathbf{C}_r$ define the  metric   $$d_r(f, g)= \sum_{n=1}^r \frac{1}{2^n} \frac{\sup_{0\le  t \le n}
	|f(t)-g(t)|}{1+\sup_{0\le  t \le n}|f(t)-g(t)|}.$$  Let $\mathfrak{C}_{\infty}$\ denote the Borel $\sigma$-field of $(\textbf{C}_{\infty},\ d)$ for $f, g \in \textbf{C}_\infty$. Note that $(\textbf{C}_{\infty},\ d)$ is a complete and separable space and the metrics $d(f, g)$ and $ d_r(f, g)$ are bounded by  1. Define mapping $Q_r $ as $Q_rf(t) = f(t), 0 \le t \le r $ for $r \ge 1$ and $f \in \textbf{C}_\infty$. We note each $Q_r$\ is $\mathfrak{C}_{\infty}\setminus \mathfrak{C}_r \ $\ measurable. Let $\mathfrak{C}_{\infty}^*$\ denote the smallest $\sigma$-field in $\textbf{C}_{\infty}$\ \textit{wrt} which  $Q_r,\ r = 1,\ 2,\ldots$\ are  measurable.
\begin{thm}
	\label{T3.1}
	(i) Fix $f \in \text{C}_{\infty}$. Then $A = B$\ where
	\begin{equation}
		\label{E3.1}
		\begin{split}
			A = \cup_{r = 1}^{\infty}\{g\ :\ g \in \textbf{C}_{\infty}.\ d_r(f,\ g) > \lambda\}\ \text{and}\\
			B = \{g\ :\ g \in \textbf{C}_{\infty},\ d(f,\ g) > \lambda\}
		\end{split}
	\end{equation}
	(ii) $\mathfrak{C}_{\infty}^* = \mathfrak{C}_{\infty}$. 
\end{thm}
\noindent	\textbf{Proof.}\\
(i)  Let $g \in B$. If it is not admitted that $g \in A$,  then $d_r(f,\ g) \le \lambda$\ for each $r \ge 1$. Since $d_r(f,\ g) \uparrow d(f,\ g)$, it follows that $d(f,\ g) \le \lambda$, a contradiction to the assumption $d(f,\ g) > \lambda$.\\
If now $g \in A$, then for some $r \ge 1$\ (and hence for all large $r$) $d_r(f,\ g) > \lambda$. Since $d(f,\ g) \ge d_r(f,\ g) > \lambda$, it follows that $g \in A$.\\
(ii) That the $Q_r$s are continuous maps is easy to verify. Hence we conclude $\mathfrak{C}_{\infty}^* \subset \mathfrak{C}_{\infty}$.\\
The reverse inclusion will stand proved if we show that every  closed $d$-sphere $S(f;\ \lambda) = \{g\ :\ g \in \textbf{C}_{\infty},\ d(f,\ g) \le \lambda\}$\  belongs to $\mathfrak{C}_{\infty}^*$. Now, since $d_r(f,\ g) \uparrow d(f,\ g)$, $S(f;\ \lambda) = \{g\ :\ d_r(f,\ g) \le \lambda\ \text{for every}\ r \ge 1\} = \cap_{r = 1}^{\infty}\{g\ :\ d_r(f,\ g) \le \lambda\}$. Since $\{Q_rg\ :\ d_r(f,\ g) \le \lambda\} \in \textbf{C}_r,\ \{g\ :\ d_r(f,\ g) \le \lambda\} \in \mathfrak{C}_{\infty}^*$. Hence $S(f;\ \lambda)$, being the intersection of a countable number of such sets, belongs to $\mathfrak{C}_{\infty}^*$.
\begin{thm}	\label{T3.2}
	(i) If $K \subset \textbf{C}_{\infty}$\ is compact, then 
	\begin{equation}\label{E3.2}
		K = \cap_{r = 1}^ {\infty}Q_r^{-1}Q_rK.
	\end{equation}
	(ii) For any set $A \subset \textbf{C}_{\infty},$
	\begin{equation}\label{E3.3}
		Q_{r + 1}^{-1}Q_{r + 1}A \subset Q_r^{-1}Q_rA.
	\end{equation}
\end{thm}
\noindent	\textbf{Proof.}\\ 
\noindent(i) That $K$\ is a subset of the rightside is trivial to see. To prove the converse, set $E = \cap_{r = 1}^ {\infty}Q_r^{-1}Q_rK.$  Then  the following relations hold.\\
$E \subset Q_r^{-1}Q_r K \;\text{for every } r \ge 1 \Rightarrow Q_r E \subset Q_r K \; \text{for every } r \ge 1 \Rightarrow Q_r^{-1}Q_r E \subset K \;\text{for every } r \ge 1$.\\
Since E is compact we  have, as observed earlier,  $ E \subset \cap_r Q_r^{-1}Q_r E $ and hence the required result follows.\\
\noindent(ii) Let $f \in Q_{r + 1}^{-1}Q_{r + 1}A$. Hence $Q_{r + 1}f \in Q_{r + 1}A$. There exists then $g \in A$\ such that $Q_{r + 1}f = Q_{r + 1}g$. This implies $Q_rf = Q_rg$\ and so $f \in Q_r^{-1}Q_rA$.
\begin{thm}	\label{T3.3}
	\(\mu_r(Q_r K),\; r = 1, 2,\ldots\)  is a monotonically decreasing sequence of real numbers.
\end{thm}
\noindent	\textbf{Proof.}\\
Let $T_r$\ map $\textbf{C}_{r + 1}$\ on to $\textbf{C}_r$\ according to the following scheme. $T_r Q_{r+1}g = Q_r g$. Thus $T_{r}Q_{r + 1} = Q_r$. Recall Wiener measure $\mu_r$\ is defined on $\mathfrak{C}_{r},\ r = 1,\ 2,\ldots $. Both $\mu_{r + 1}T_r^{-1}$\ and $\mu_r$\ are measures defined on $\mathfrak{C}_r$. Their finite dimensional distributions are the same. Hence the two are identical (ref. Remark \ref{R2.3}). We then have
(using the formula for change of variables in an integral (ref. Theorem 2.3.6, p91,Pakshirajan , 2013), $\mu_r(Q_rK) = \int\limits_{Q_rK}\;d\mu_r = \int\limits_{Q_rK}\;d\mu_{r + 1}T_r^{-1} = \int\limits_{T_r^{-1}Q_rK}\;d\mu_{r + 1} \ge \int\limits_{Q_{r + 1}K}\;d\mu_{r + 1}$. We see from this that $\mu_r(Q_rK)$\ is a monotonically decreasing sequence of real numbers. \\
\noindent	For $K \in \mathfrak{C}_{\infty},\ K$\ compact, define
\begin{equation}
	\label{E3.3}
	\mu_{\infty}K = \lim\limits_{r \uparrow \infty}\mu_r(Q_rK)
\end{equation}  and for arbitrary $A \in \mathfrak{C}_{\infty}$, define
\begin{equation}	\label{E3.4}
	\mu_{\infty}A = \sup\limits_{K \subset A,\ K \text{compact}}\mu_{\infty}K
\end{equation}
and proceed as in the construction of the measure $\mu_a$, use (\ref{E3.3}) and arrive at a countably additive finite measure $\mu_{\infty}$  which is finite and $\le 1$, by(\ref{E3.3}) and (\ref{E3.4}). \\
That $\mu_{\infty}$\ is a probability measure will follow if we show that $\mu_{\infty}\textbf{C}_{\infty} = 1$. \\ 
Consider the  H$\ddot{o}$lder space $H_{\alpha, \infty}$ and define, for $x, y \in H_{\alpha, \infty}$ the metric 
\begin{equation*}
	d_{\alpha, \infty}(x, y)=\sum_{n=1}^\infty\frac{1}{2^n}\frac{\sup_{0\le t, s \le n; 0<|t-s|<1}\frac{|x(t)-y(t)-x(s)+y(s)|}{|t-s|^\alpha}}{1+\sup_{0\le t,s \le n; 0< |t-s| < 1} \frac{|x(t)-y(t)-x(s)+y(s)|}{|t-s|^\alpha}}
\end{equation*}

Also define 
\begin{equation*}
	d^*_{\alpha, \infty}(x, y)=\sum_{n=1}^\infty\frac{1}{2^n}\frac{\sup_{0\le t, s \le n; 0<|t-s|< 1} \frac{|x(t)-y(t)|}{|t-s|^\alpha}}{1+\sup_{0\le t, s \le n;0<|t-s|< 1}\frac{|x(t)-y(t)|}{|t-s|^\alpha}}.
\end{equation*}
Note that $	d^*_{\alpha, \infty}(x, y) \le 	d_{\alpha, \infty}(x, y) \le 2 	d^*_{\alpha, \infty}(x, y).$
Further define on $H_{\alpha, \infty}$ another metric 
\begin{equation*}
	d(x, y)= \sum_{n=1}^\infty\frac{1}{2^n}\frac{\sup_{0 \le t \le n}|x(t)-y(t)|}{1+\sup_{0 \le t \le n}|x(t)-y(t)|}.
\end{equation*}
Note that $d(x, y) \le d^*_{\alpha, \infty}(x, y) \le 	d_{\alpha, \infty}(x, y)$
\noindent Denote by $\vartheta$\ the null element. i.e., the function that is identically zero. Since $S^*(\lambda) = \{f\ :\ d(f,\ \vartheta) \le \lambda\} \uparrow \textbf{C}_{\infty}\ \text{as}\ \lambda \uparrow \infty$, it is sufficient to show that, given $\varepsilon > 0$, a 
$\lambda$\ can be found such that $\mu_{\infty}\big(S^*(\lambda)\big) > 1 - \varepsilon$. \\ 
Define $ \mathbb{S}^*_{\alpha}(\lambda)=\{x: x\in H_{\alpha, \infty}, d_{\alpha, \infty}(\vartheta, x) \le \lambda\}.  $	Since $d(x,\ y) \le d_{\alpha,\infty}(x, y)$, $\mathbb{S}^*_{\alpha}(\lambda) \subset S^*(\lambda)$. Hence it is enough to find a $\lambda$\ such that $\mu_{\infty}(\mathbb{S}^*_{\alpha}(\lambda)) > 1 - \varepsilon$. Since $\mathbb{S}^*_{\alpha}(\lambda)$\ is a compact set, $\mu_{\infty}(\mathbb{S}^*_{\alpha}(\lambda)) = \lim\limits_{r  \rightarrow \infty}\mu_r(Q_r\mathbb{S}^*_{\alpha}(\lambda))$ by(\ref{E3.3}). Take $r=[\lambda].$ The arguments in the proof of Theorem \ref{T2.3} apply and we get $\mu_{\infty}(\mathbb{S}^*_{\alpha}(\lambda)) \ge e^{-r/\lambda^2}\ge e^{-1/\lambda} \rightarrow 1$\ as $\lambda \rightarrow \infty$.
\section{appendix}
Recall
 \begin{equation*}
	\mathbb{S}_{\alpha, a}(\lambda)=\{ x: \|x\|_\alpha=\sup_{0< s, t \le a;\; 0 < |s-t| < 1}\frac{|x(t)-x(s)|}{|s-t|^\alpha} \le \lambda\}.
\end{equation*}
We shall show that  $	\mathbb{S}_{\alpha, a}(\lambda)$ is a compact subset of $\bf{C}_a$. To this end we shall show that (a) that it is bounded, (b) that it is closed and (c) that it is uniformly
equicontinuous.\\
Let $\alpha < \beta.$ Then $\mathbb{S}_{\beta, a}(\lambda) \subset \mathbb{S}_{\alpha, a}(\lambda)$.\\\\
(b) Consider $ u \in H^0_\alpha,\; u_n \in \mathbb{S}_{\beta, a}(\lambda)$ and let $\|u_n -u\|_\alpha \rightarrow 0$ as $n \rightarrow \infty$. 
\begin{equation*}
	\begin{split}
		\sup_{0< s, t \le a;\; 0 < |s-t| < 1}\frac{|u(t)-u(s)|}{|s-t|^\beta} \le \\
		\le \bigg(\sup_{0< s, t \le a;\; 0 < |s-t| < 1}\frac{|u(t)-u(s)|}{|s-t|^\alpha}\bigg)^{\beta/\alpha}\bigg(\sup_{0< s, t \le a;\; 0 < |s-t| < 1}|u(t)-u(s)| \bigg)^{1-\frac{\beta}{\alpha}}\\
		=\lim_{n\rightarrow \infty}\bigg(\sup_{0< s, t \le a;\; 0 < |s-t| < 1}\frac{|u_n(t)-u_n(s)|}{|s-t|^\alpha}\bigg)^{\beta/\alpha}\bigg(\sup_{0< s, t \le a;\; 0 < |s-t| < 1}|u(t)-u(s)| \bigg)^{1-\frac{\beta}{\alpha}}\\ 
		\le \limsup_{n\rightarrow \infty}\bigg(\sup_{0< s, t \le a;\; 0 < |s-t| < 1} \frac{|u_n(t)-u_n(s)|}{|s-t|^\beta}\bigg)\times \\\limsup_{n\rightarrow \infty} \sup_{0< s, t \le a;\; 0 < |s-t| < 1} \bigg( |u_n(t)-u_n(s)|^{\frac{\beta}{\alpha}-1}|u(t)-u(s)|^{1-\frac{\beta}{\alpha}}\bigg)\\
		\le \limsup_{n \rightarrow \infty}\|u_n\|_\beta \le \lambda.
	\end{split}
\end{equation*}
This shows that $u\in \mathbb{S}_{\beta, a}.$ i.e.,  $\mathbb{S}_{\beta, a}(\lambda) $ is a closed subset of $H^0_\alpha$.\\
(c) Finally we prove that $\mathbb{S}_{\alpha, a}$ is uniformly equicontinuous in $H^0_{\alpha}$. Let  $x \in \mathbb{S}_{\alpha, a}(\lambda)$ and let $\delta > 0$ be arbitrary. Consider , for $\alpha < \beta$
\begin{equation*}
	\sup_{0 \le s, t \le a;\; |t-s| \le \delta}\frac{|x(t)-x(s)|}{|t-s|^\alpha} = \sup_{0 \le s, t \le a;\; |t-s| \le \delta}\frac{|x(t)-x(s)|}{|t-s|^\beta}\times|t-s|^{\beta -\alpha} \le \lambda \;\delta^{\beta-\alpha}.
\end{equation*}
Since $\delta $ is arbitrary and since the above step holds for all $x \in \mathbb{S}_{\alpha, a}$ uniform equicontinuity  follows. \\
We thus have compactness of $\mathbb{S}_{\alpha, a}$.

\end{document}